\documentclass [12pt]{article}
\usepackage{amssymb}
\def\lanbox{\hbox{$\, \vrule height 0.25cm width 0.25cm depth 0.01cm \,$}}

\usepackage{color}

\begin{document}

\begin{center}
{\Large\bf Stability and instability  of constant stationary  
\\
solutions to some integro-differential equations }
\end {center}

\bigskip 

\centerline{Yuming Chen$^{1}$, Vitali Vougalter$^{2 \ *}$}

\vspace*{0.5cm}

\centerline{$^1$ Department of Mathematics, Wilfrid Laurier University}

\centerline{Waterloo, Ontario, N2L 3C5, Canada}

\centerline{ e-mail: ychen@wlu.ca}

\medskip

\centerline{$^{2 \ *}$ Department of Mathematics, University of Toronto}

\centerline{Toronto, Ontario, M5S 2E4, Canada}

\centerline{e-mail: vitali@math.toronto.edu}

\bigskip 

\noindent {\bf Abstract.}  We consider some reaction-diffusion equations describing
systems with the nonlocal consumption of resources and the intraspecific
competition. Sharp conditions on the coefficients are obtained to  ensure the stability and instability of nontrivial constant stationary solutions.

\medskip



\bigskip
\bigskip
\bigskip
\bigskip

\setcounter{equation}{0}

\section{Introduction}

In the first part of this article we consider the following integro-differential equation,
\begin{equation}
\label{id1}
\frac{\partial u}{\partial t}=\Delta u+ku^{2}\left(1-au-b\int_{{\mathbb R}^{n}}
\phi(x-y)u(y,t)dy \right),
\end{equation}
where
$(x,t)\in {\mathbb R}^{n}\times {\mathbb R}$, $n\in {\mathbb N}$, $u\in
\mathcal{C}:=C(\mathbb{R}^n,\mathbb{R})$. Here the
integral kernel $\phi(x)\geq 0$ for $x\in {\mathbb R}^{n}$, which is
dependent only upon the radial variable $|x|$ or it is a product of even
functions (see e.g. expressions~(\ref{fi51}), (\ref{fi53}), (\ref{fi54}))
and is normalized as $\int_{{\mathbb R}^{n}}\phi(x)dx=1$. The parameters
$k$, $b>0$ and $a\geq 0$. Obviously, the only nontrivial solution constant
in space and time the equation~(\ref{id1})  possesses is
\begin{equation}
\label{ss0}
u(x,t)=\frac{1}{a+b}.
\end{equation}
Linearization of problem~(\ref{id1}) around~(\ref{ss0})  with the ansatz
$u(x,t)=\frac{1}{a+b}+\psi(x)e^{-{\lambda t}}$ yields the
spectral problem on $L^{2}({\mathbb R}^{n})$
\begin{equation}
\label{sppr0}
{\cal L}\psi=\lambda \psi,
\end{equation}
where
\[
{\cal L}\psi=-\Delta \psi+k_{1}\psi+k_{2}\int_{{\mathbb R}^{n}}\phi(x-y)
\psi(y)dy
\]
with $k_{1}=\frac{ka}{(a+b)^{2}}$ and $
k_{2}=\frac{kb}{(a+b)^{2}}$. The spectrum of the operator $\cal L$ is  denoted as $\sigma({\cal L})$. Problem~(\ref{sppr0}) is said to be {\it spectrally stable} if
$\sigma({\cal L})\subseteq [0, \infty)$ and {\it spectrally unstable}
if there exists $\lambda\in\sigma({\cal L})$  such that $\lambda<0$.

We denote the standard Fourier transform of a function
$\psi(x)\in L^{2}({\mathbb R}^{n})$ as
\[
\widehat{\psi}(p)=\frac{1}{ (2\pi)^{n\over 2}}\int_{{\mathbb R}^{n}}
\psi(x)e^{-ip\cdot x}dx, \quad p\in {\mathbb R}^{n},
\]
which enables us to write the quadratic form of the linearized operator as
\[
({\cal L}\psi, \psi)=\int_{{\mathbb R}^{n}}\Phi(p)|\widehat{\psi}(p)|^{2}dp,
\]
where
\[
\Phi(p)=|p|^{2}+k_{1}+k_{2}(2\pi)^{\frac{n}{2}}\widehat{\phi}(p).
\]
Clearly problem~(\ref{sppr0}) is spectrally stable if for the Fourier image
of the integral kernel involved in equation~(\ref{id1}) we have
$\widehat{\phi}(p)\geq 0$ for $p\in {\mathbb R}^{n} \ \mathrm{a.e.}$.
Spectral instability arises when there is a region
$\Omega\subset {\mathbb R}^{n}$ with nonempty interior
such that $\Phi(p)<0$ for $p\in \Omega \ \mathrm{a.e.}$, which
occurs  due to the negativity of $\widehat{\phi}(p)$. Then we can choose a
trial function $\psi(x)\in L^{2}({\mathbb R}^{n})$ such that its Fourier
image $\widehat{\psi}(p)$ is
compactly supported inside $\Omega$. The negativity of the resulting quadratic
form $({\cal L}\psi, \psi)$ along with the min-max principle
(see e.g.\ \cite[p.306]{LSW76})
implies that the operator ${\cal L}$ has negative spectrum. 

We consider the following types of integral kernels whose Fourier
transforms can be easily computed. The first one is  
\begin{equation}
\label{fi1}
\phi_{1}(x)=\frac{\alpha}{2}e^{-\alpha|x|}, \qquad x\in {\mathbb R},  \alpha>0 
\end{equation}
and hence  $\widehat{\phi_{1}}(p)=\frac{\alpha^{2}}{\sqrt{2\pi}(p^{2}+
{\alpha}^{2})}$, $p\in {\mathbb R}$; the second one is the  two
dimensional generalization of $\phi_1$,
\begin{equation}
\label{fi2}
\phi_{2}(x,y)=\phi_{1}(x)\phi_{1}(y)=\frac{\alpha^{2}}{4}e^{-\alpha|x|}
e^{-\alpha|y|}, \qquad x,y \in {\mathbb R},  \alpha>0
\end{equation}
with  $\widehat{\phi_2}(p)=\frac{\alpha^4}{2\pi(p_1^2+\alpha^2)(p_2^2+\alpha)}$, $p=(p_1,p_2)\in\mathbb{R}^2$; the third one is the Gaussian kernel
\begin{equation}
\label{fi3}
\phi_{3}(x)=\left(\frac{\alpha}{\pi}\right)^{\frac{n}{2}}e^{-\alpha |x|^{2}},
\qquad x\in {\mathbb R}^{n},  n\in {\mathbb N}, \alpha>0.
\end{equation}
One can get $\widehat{\phi_{3}}(p)=
\frac{1}{(2\pi)^{\frac{n}{2}}}e^{-\frac{|p|^{2}}{4\alpha}}, p\in {\mathbb R}^{n}$. The
fourth one is the three dimensional generalization of $\phi_1$, \begin{equation}
\label{fi4}
\phi_{4}(x)=\frac{\alpha^{3}}{8\pi}e^{-\alpha |x|}, \qquad x\in {\mathbb R}^{3}, \alpha>0.
\end{equation}
Thus using the formula for the heat kernel of the square root of the
Laplacian (see e.g.\ \cite[p.169]{LL97}) we easily obtain
$ \widehat{\phi_{4}}(p)=
\frac{\alpha^{4}}{(2\pi)^{\frac{3}{2}}({\alpha}^{2}+|p|^{2})^{2}}$, $p\in {\mathbb R}^{3}$. Hence  these choices of
kernels~(\ref{fi1})-(\ref{fi4}) having positive Fourier
images in the whole space yield spectral stability for problem~(\ref{sppr0}).
On the other hand if we choose the integral kernel
\begin{equation}
\label{fi5}
\phi_{5}(x)=\frac{1}{2N}\chi_{[-N,\ N]}(x), \qquad x\in {\mathbb R}, \quad
N>0,
\end{equation}
where $\chi_{A}$ stands for the characteristic function of a set $A$, we
arrive at
\[
\widehat{\phi_{5}}(p)=\frac{1}{\sqrt{2 \pi}}\frac{\sin pN}{pN},
\qquad p\in {\mathbb R},
\]
which is sign indefinite. Now we form products out of
the kernels above, namely,
\begin{eqnarray}
\label{fi51}
\phi_{5}(x)\phi_{1}(y) & =& \frac{1}{2N}\chi_{[-N,\ N]}(x)
{\alpha\over 2}e^{-\alpha |y|}, \qquad x,y \in {\mathbb R},
\\
\label{fi53}
\phi_{5}(x)\phi_{3}(y)& = & \frac{1}{2N}\chi_{[-N, \ N]}(x)
\sqrt{\frac{\alpha}{\pi}}e^{-\alpha y^{2}}, \qquad x,y \in {\mathbb R},
\\
\label{fi54}
\phi_{5}(x)\phi_{4}(y) & = & \frac{1}{2N}\chi_{[-N, \ N]}(x)
\frac{\alpha^{3}}{8 \pi}e^{-\alpha |y|}, \qquad x\in {\mathbb R}, y\in {\mathbb R}^{3}.
\end{eqnarray}
Since the linearized operator involved in~(\ref{sppr0}) is similar to
the one studied in~\cite{ABVV10} (see also~\cite{BNPR09,GVA06}),
we arrive at the following statement analogous to Theorem 2.8 of~\cite{ABVV10} (see also Chapter 3
of~\cite{GVA06}).

\medskip

\noindent {\bf Theorem 1.} {\it Let $k_{1}=0$ and the kernel $\phi$ be one of~(\ref{fi5}), (\ref{fi51}), (\ref{fi53}) and (\ref{fi54}). Then
problem~(\ref{sppr0}) is spectrally stable if
$\frac{1}{k_{2}}\geq -N^{2}{\frac{\sin z_{1}}{z_{1}^{3}}}$
and is spectrally unstable if $\frac{1}{k_{2}}< -N^{2}\frac{\sin z_{1}}{z_{1}^{3}}$, where
$z_{1}$ is the unique solution of the equation $\tan z=\frac{1}{3}z$ on
the interval $(\pi, \frac{3\pi}{2})$.}
 
Emergence and propagation of patterns in nonlocal reaction-diffusion
equations arising in the theory of speciation was covered in~\cite{VV13}.
Doubly nonlocal reaction-diffusion equations and the emergence of species
were  discussed in~\cite{BVV17}. Asymptotic analysis of SIR epidemic model with nonlocal diffusion and general  nonlinear incidence was performed
in~\cite{DCB23}. The article~\cite{KLS23} deals with the
diffusion-driven blow-up for a nonlocal Fisher-KPP type model. A
discretization method for nonlocal diffusion type equations was developed in
~\cite{MORV22}.

In the second part of this article we study the following integro-differential equation,
\begin{equation}
\label{id2}
\frac{\partial u}{\partial t}=d\Delta u+u^{2}\left(a
-\int_{{\mathbb R}^{n}}
\phi(x-y)u(y,t)dy \right)-bu,
\end{equation}
where $(x,t)\in {\mathbb R}^{n}\times {\mathbb R}$,  $n\in {\mathbb N}$, $u\in
\mathcal{C}$, the integral kernel $\phi$ satisfies the same requirements as
those  for~(\ref{id1}), the diffusion coefficient
$d>0$ and the parameters $a$, $b>0$. We assume that
$\frac{a^{2}}{4}\geq b$ such that equation~(\ref{id2})
admits two nontrivial real-valued solutions constant in space and time
given by
\begin{equation}
\label{ss12}
c_{1}=\frac{a}{2}+\sqrt{\frac{a^{2}}{4}-b}, \quad
c_{2}=\frac{a}{ 2}-\sqrt{\frac{a^{2}}{4}-b}
\end{equation}
and, in particular $c_{1,2}=\frac{a}{2} $ when
$\frac{a^{2}}{4}=b$.
Linearization of equation~(\ref{id2}) around a solution in~(\ref{ss12}) with the ansatz
$u(x,t)=c_{1,2}+\psi(x)e^{-{\lambda t}}$ yields the
spectral problem on $L^{2}({\mathbb R}^{n})$
\begin{equation}
\label{spprk}
{\cal L}_{k}\psi=\lambda \psi, \quad k=1,2,
\end{equation}
where
\[
{\cal L}_{k}\psi=-d \Delta \psi+c_{k}^{2}\int_{{\mathbb R}^{n}}\phi(x-y)\psi(y)dy-
b\psi, \quad k=1,2,
\]
such that problem~(\ref{spprk}) is spectrally stable if
$\sigma({\cal L}_{k})\subseteq[0, \infty)$ and is spectrally unstable if there
exists $\lambda\in \sigma({\cal L}_{k})$ such that $\lambda<0$. Using the Fourier
transform we write the quadratic form
\[
({\cal L}_{k}\psi, \psi)=\int_{{\mathbb R}^{n}}\Phi_{k}(p)|\widehat{\psi}(p)|^{2}dp,
\]
where
\[
\Phi_{k}(p)=d|p|^{2}+(2 \pi)^{\frac{n}{2}}c_{k}^{2}\widehat{\phi}(p)-b, \quad
k=1,2.
\]
By considering the integral kernels mentioned above we show that problem~(\ref{spprk}) is distinct from~(\ref{sppr0}). In fact, (\ref{spprk}) can be spectrally unstable
even if the Fourier image of the integral kernel involved  in~(\ref{id2}) is positive in the whole space. More precisely, we have the following results.

\medskip 

\noindent {\bf Theorem 2.} {\it Let the kernel $\phi$ be given by (\ref{fi1}).}
\begin {itemize}
\item [{\rm (i)}] \textit {When $\frac{a^{2}}{4}>b$, problem~(\ref{spprk}) with
$k=1$ is spectrally stable if  $ d\geq \frac{(c_{1}-\sqrt{c_{1}^{2}-b})
^{2}}{\alpha^{2}}$
and is spectrally unstable if
$0<d<\frac{(c_{1}-\sqrt{c_{1}^{2}-b})^{2}}{\alpha^{2}}$, and  problem (\ref{spprk}) with $k=2$ is spectrally unstable.}

\item [{\rm (ii)}] \textit {When $\frac{a^{2}}{4}=b$ problem (\ref{spprk}) with
$k=1$, $2$ is spectrally stable if
$d\geq \frac{b}{\alpha^{2}}$ and is spectrally unstable if
$0< d<\frac{b}{\alpha^{2}}$.}
\end {itemize}

\noindent {\bf Theorem 3.} {\it Let the kernel $\phi$ be given by (\ref{fi2}).

\begin {itemize}
\item [{\rm (i)}] When $\frac{a^2}{4}>b$, problem (\ref{spprk}) with
$k=1$ is spectrally stable if
$ d\geq \frac{x^{*}}{\alpha^{2}}$
and is spectrally unstable if
$0<d<\frac{x^{*}}{\alpha^{2}}$, where $x^{*}$ is defined
in Lemma~7, and 
Problem (\ref{spprk}) with $k=2$ is spectrally unstable.

\item [{\rm (ii)}] When $\frac{a^{2}}{ 4}=b$, problem (\ref{spprk}) with
$k=1$, $2$ is spectrally stable if
$d\geq \frac{a^{2}}{4\alpha^{2}}$ and is spectrally unstable if
$0<d< \frac{a^{2}}{4\alpha^{2}}$.
\end {itemize}
}

\noindent {\bf Theorem 4.} {\it Let the kernel $\phi$ be given by (\ref{fi3}).
\begin{itemize}
\item  [{\rm (i)}] When $\frac{a^{2}}{4}>b$, problem (\ref{spprk}) with
$k=1$ is spectrally stable if
$ d\geq s_{0}(a,b)\frac{c_{1}^{2}}{4\alpha}$
and is spectrally unstable if
$0<d< s_{0}(a,b)\frac{c_{1}^{2}}{4\alpha}$, where $s_{0}(a,b)$
is defined in Lemma~8, and problem (\ref{spprk}) with $k=2$ is spectrally unstable.

\item [{\rm (ii)}] When $\frac{a^{2}}{4}=b$ problem (\ref{spprk}) with
$k=1$, $2$ is spectrally stable if
$d\geq \frac{b}{4\alpha}$ and is spectrally unstable if
$0<d< \frac{b}{ 4\alpha}$.
\end {itemize}}

\noindent {\bf Theorem 5.} {\it Let the kernel $\phi$ be given by (\ref{fi4}).
\begin {itemize}
\item [{\rm (i)}] When $\frac{a^{2}}{4}>b$, problem (\ref{spprk}) with
$k=1$ is spectrally stable if
$d\geq \frac{2x^{*}}{\alpha^{2}}$
and is spectrally unstable if
$0<d<\frac{2x^{*}}{\alpha^{2}}$, where $x^{*}$ is defined
in Lemma 7, and  problem (\ref{spprk}) with $k=2$ is spectrally unstable.

\item [{\rm (ii)}] When $\frac{a^{2}}{4}=b$, problem (\ref{spprk}) with
$k=1$, $2$ is spectrally stable if
$d\geq \frac{a^{2}}{2\alpha^{2}}$ and is spectrally unstable if
$0<d< \frac{a^{2}}{2\alpha^{2}}$.
\end {itemize}
}

These results are proved in the next section.


\setcounter{equation}{0}

\section{Proofs of Theorems 2--5}

For technical purpose, we first establish the following a little bit  trivial result, which
proves the statements of Theorems 2--5 in the case where
$\frac{a^2}{4}>b$ and $k=2$.

\medskip

\noindent {\bf Lemma 6.} {\it When $\frac{a^{2}}{4}>b$, we have
$c_{1}^{2}-b>0$ and $ c_{2}^{2}-b<0$ and hence problem (\ref{spprk}) with
the kernels given by (\ref{fi1})--(\ref{fi4}) and $k=2$ is spectrally
unstable. Moreover,  when $\frac{a^{2}}{4}=b$ both $c_{1}^{2}-b$ and $c_{2}^{2}-b$ are equal to zero.
}

{\it Proof.} By means of (\ref{ss12}), 
\[
c_{1}^{2}-b=\left(\frac{a}{2}-\sqrt{b}+\sqrt{\frac{a^{2}}{4}-b}\right)
\left(c_{1}+\sqrt{b}\right),
\]
which is positive if  $\frac{a^{2}}{4}>b$ and vanishes if
$\frac{a^{2}}{4}=b$. Similarly
\[
c_{2}^{2}-b=\left(\frac{a}{2}-\sqrt{b}-\sqrt{\frac{a^{2}}{4}-b}\right)
\left(\frac{a}{2}+\sqrt{b}-\sqrt{\frac{a^{2}}{4}-b}\right),
\]
which clearly vanishes when $\frac{a^{2}}{4}=b$.
For $\frac{a^{2}}{4}>b$, we expand the second term in the
product above as
\[
\sqrt{\frac{a}{2}+\sqrt{b}}\left(\sqrt{\frac{a}{2}+\sqrt{b}}-
\sqrt{\frac{a}{2}-\sqrt{b}}\right)>0
\]
and the first one as
\[
\sqrt{\frac{a}{2}-\sqrt{b}}\left(\sqrt{\frac{a}{2}-\sqrt{b}}-
\sqrt{\frac{a}{2}+\sqrt{b}}\right)<0.
\]
Since $\int_{{\mathbb R}^{n}}\phi(x)dx=1$, we have
$\widehat{\phi}(0)=\frac{1}{(2\pi)^{n\over 2}}$ and it follows that $\Phi_{k}(0)=c_{2}^{2}-b<0$ when $k=2$ and $\frac{a^{2}}{4}>b$. Hence by continuity for the kernels given by
formulas (\ref{fi1})--(\ref{fi4}), the function $\Phi_{2}(p)$ is negative
in a neighborhood $\Omega$ of the origin with nonempty interior. Let us
choose the trial function $\widehat{\psi}(p)\in L^{2}({\mathbb R}^{n})$ to be compactly supported inside $\Omega$. Then the quadratic form
$({\cal L}_{2}\psi, \psi)<0$ which by means of the min-max principle implies the existence of negative spectrum and spectral instability of the problem~(\ref{spprk}) with $k=2$ and the kernels given by (\ref{fi1})--(\ref{fi4}).
$\hfill\lanbox$

\smallskip

We first prove Theorem 2.

{\it Proof of Theorem 2.} When the kernel is given by (\ref{fi1}), we have
\[
\Phi_{k}(p)=dp^{2}+c_{k}^{2}\frac{\alpha^{2}}{p^{2}+\alpha^{2}}-b, \qquad
p\in {\mathbb R},  k=1,2.
\]

Let us first consider the case of two distinct constant nontrivial stationary solutions of equation~(\ref{id2}), ie., $\frac{a^2}{4}>b$. When $k=1$, clearly, it is sufficient to  investigate the existence of
intervals of negativity for the function
\[
N_{1}(p)=dp^{4}+(d\alpha^{2}-b)p^{2}+\alpha^{2}(c_{1}^{2}-b), \qquad
p\in {\mathbb R}.
\]
We distinguish the two cases of large and small  diffusion
coefficient. If $d\geq \frac{b}{\alpha^{2}}$, by Lemma 6, $N_1(p)$  is strictly positive on ${\mathbb R}$. This implies
spectral stability for problem (\ref{spprk}) with $k=1$. Now suppose $0<d< \frac{b}{\alpha^{2}}$. Then $N_{1}(p)$
is minimal at $p^{2}=-\frac{d\alpha^{2}-b}{2d} $ and a
straightforward computation yields the minimal value
\[
\frac{\alpha^{2}}{4d}\left[-\alpha^{2}d^{2}+4d\left (c_{1}^{2}-\frac {b}{2}\right)-
\frac{b^{2}}{\alpha^{2}}\right].
\]
The quadratic polynomial in $d$ inside the brackets above has two roots
\[
d_{1}=\frac{(c_{1}+ \sqrt{c_{1}^{2}-b})^{2}}{\alpha^{2}}, \quad
d_{2}=\frac{(c_{1}- \sqrt{c_{1}^{2}-b})^{2}}{\alpha^{2}}.
\]
It can be easily verified that
$0<d_{2}<\frac{b}{\alpha^{2}}<d_{1}$. Therefore, when
$d\geq d_{2}$ we have $\Phi_{1}(p)$ nonnegative on ${\mathbb R}$ and therefore,
spectral stability for problem (\ref{spprk}) with $k=1$. When $0<d<d_{2}$
there are intervals of negativity of $\Phi_{1}(p)$ on the real line, which
by means of the min-max principle yields the existence of negative spectrum
and therefore, spectral instability for problem (\ref{spprk}) with $k=1$.

Finally, we study the situation when both nontrivial constant solutions
of (\ref{id2}) coincide, i.e., $\frac{a^2}{4}=b$. Then  $c_{1}=c_{2}=\frac{a}{2}$. In this case
\[
\Phi_{k}(p)=\frac{dp^{2}+d{\alpha}^{2}-b}{1+\frac{{\alpha}^{2}}{ p^{2}}}, \qquad
p\in {\mathbb R}, k=1,2.
\]
Hence for $d\geq \frac{b}{{\alpha}^{2}}$, $\Phi_k(p)$  is
nonnegative on ${\mathbb R}$ which implies spectral stability for our problem
and for $0<d< \frac{b}{{\alpha}^{2}}$ it is negative in the
neighborhood of the origin, which means spectral instability for (\ref{spprk}).
$\hfill\lanbox$

\smallskip

Before moving on, we need the following technical lemma.

\noindent {\bf Lemma 7.} {\it Assume $\frac{a^2}{4}>b$. Consider the function $
f(x)=27c_{1}^{2}x^{2}-(b+2x)^{3}$, $x\in (0, b)$. $f$ has the
minimum only  at
$x_{1}=\frac{9c_{1}^{2}-4b-3c_{1}\sqrt{9c_{1}^{2}-8b}}{8}$
and is monotonically increasing on $(x_{1},b)$.
Moreover, in the interval $(0,b)$ the equation $f(x)=0$  admits the
unique solution $x^{*}\in (x_{1}, b)$.}

{\it Proof.} A simple calculation yields
\[
f'(x)=-6[4x^{2}-(9c_{1}^{2}-4b)x+b^{2}],
\]
which has two real zeros via Lemma 6
\[
x_{1}=\frac{9c_{1}^{2}-4b-3c_{1}\sqrt{9c_{1}^{2}-8b}}{8}, \quad
x_{2}=\frac{9c_{1}^{2}-4b+3c_{1}\sqrt{9c_{1}^{2}-8b}}{ 8}.
\]
Since both the product and the sum of $x_{1}$ and $x_{2}$ are positive,
we have $x_{1, 2}>0$. Apparently, it can be verified that
$0<x_{1}<b<x_{2}$. Near the endpoints of the $(0, b)$ interval
$f(0^{+})=-b^{3}<0$ and $f(b^{-})=27b^{2}(c_{1}^{2}-b)>0$ via Lemma 6.
Thus $f(x)$ considered on $(0, \ b)$ attains its minimal value, which is
negative at $x_{1}$, is strictly increasing for $x>x_{1}$, and has the
unique zero at some point $x^{*}\in (x_{1}, b)$.
$\hfill\lanbox$

\medskip

Next, we prove Theorem 3.

{\it Proof of Theorem 3.} When the kernel function is given by (\ref{fi2}),
we have
\begin{equation}
\label{Fik3}
\Phi_{k}(p)=d|p|^{2}+c_{k}^{2}\frac{\alpha^{4}}{ (p_{1}^{2}+{\alpha}^{2})
(p_{2}^{2}+{\alpha}^{2})}-b, \quad p=(p_{1}, p_{2})\in {\mathbb R}^{2},
k=1,2.
\end{equation}

We start  with the case of
$\frac{a^{2}}{4}>b$ and therefore, problem (\ref{id2}) has two  distinct constant  stationary
solutions $c_{1}$ and $c_{2}$.
When $k=1$, in  polar coordinates $p=(|p|, \theta)$,  $|p|\in [0,  \infty)$, $\theta\in [0, 2\pi)$, the function $\Phi_{1}(p)$
obviously equals  
\[
d|p|^{2}+c_{1}^{2}\frac{\alpha^{4}}{\frac{1}{ 4}|p|^{4}\sin^{2}2\theta+
{\alpha}^{2}|p|^{2}+\alpha^{4}}-b,
\]
which is minimal for
$\theta=\frac{\pi}{4}$, $\frac{3\pi}{4}$, $\frac{5\pi}{4}$, $\frac{7\pi}{4}$, namely when $p_{1}^{2}=p_{2}^{2}$. In this case we introduce
the new variable $q=p_{1}^{2}\in [0, \infty)$ and arrive at the expression
\[
\frac{(2dq-b)(q+{\alpha}^{2})^{2}+c_{1}^{2}{\alpha}^{4}}{ (q+{\alpha}^{2})^{2}}.
\]
Since its denominator is strictly positive, it makes sense to consider
\[
N_{2}(q)=(2dq-b)(q+{\alpha}^{2})^{2}+c_{1}^{2}{\alpha}^{4},
\]
which is positive at the origin by means of Lemma 6.
We easily get
\[
\frac{dN_{2}(q)}{dq}=2(q+{\alpha}^{2})(3dq+d{\alpha}^{2}-b).
\]
Thus when the diffusion coefficient is large enough, namely
$d\geq \frac{b}{{\alpha}^{2}}$, the function $N_{2}(q)$ is
monotonically increasing on the interval $(0,\infty)$. This implies its
strict positivity and the spectral stability of $(\ref{spprk})$ with $k=1$.
For  small values of the diffusion coefficient,
$0<d< \frac{b}{{\alpha}^{2}}$ the function $N_{2}(q)$ attains
its minimal value at
\[
q_{0}=\frac{b-d\alpha^{2}}{3d}.
\]
A straightforward computation yields
\[
N_{2}(q_{0})=\frac{27c_{1}^{2}{\alpha}^{4}d^{2}-(b+2d{\alpha}^{2})^{3}}
{27d^{2}}.
\]
Let us introduce the new variable $x=d\alpha^{2}$, $0<x<b$. This will turn
the numerator of the expression in the right side of $N_2(q_0)$  
into the  function $f(x)$ studied in Lemma 7. Therefore, problem (\ref{spprk})
with $\frac{a^{2}}{4}>b$ and $k=1$ is spectrally stable for
$d\geq \frac{x^{*}}{{\alpha}^{2}}$ and is spectrally unstable
when $0<d< \frac{x^{*}}{{\alpha}^{2}}$.

Now consider the  degenerate case that $\frac{a^{2}}{4}=b$. In this case, 
\begin{equation}
\label{Fi123}
\Phi_k(p)=d|p|^{2}+\frac {a^{2}}{4}\frac{{\alpha}^{4}}{ (p_{1}^{2}+{\alpha}^{2})
(p_{2}^{2}+{\alpha}^{2})}-\frac{a^{2}}{4}, \quad p=(p_{1},p_{2})\in
{\mathbb R}^{2}, k=1,2.
\end{equation}
By changing the variables to the polar coordinates we easily obtain
\[
\Phi_k(p)=d|p|^{2}+\frac{a^{2}{\alpha}^{4}}{ |p|^{4}sin^{2}2\theta+4{\alpha}^{2}|p|^{2}+
4{\alpha}^{4}}-\frac{a^{2}}{4},
\]
which is minimal when $p_{1}^{2}=p_{2}^{2}$. Using the variable
$q=p_{1}^{2}\in [0, \ \infty)$ we arrive at
\[
\frac{(8dq-a^{2})(q+{\alpha}^{2})^{2}+a^{2}{\alpha}^{4}}{ 4(q+{\alpha}^{2})^{2}}.
\]
Due to the positivity of the denominator in the expression above, one needs
to study
\[
N_{3}(q)=(8dq-a^{2})(q+{\alpha}^{2})^{2}+a^{2}{\alpha}^{4}, \quad q\in
[0,  \infty).
\]
Obviously, $N_{3}(0)=0$. A straightforward computation yields
\[
\frac{dN_{3}(q)}{dq}=2(q+{\alpha}^{2})(12dq+4d{\alpha}^{2}-a^{2}).
\]
Therefore, when the diffusion coefficient
$d\geq \frac{a^{2}}{4{\alpha}^{2}}$, we have
$N_{3}(q)\geq N_{3}(0)=0$, $q\in [0, \infty)$ and hence problem
(\ref{spprk}) with $\frac{a^{2}}{4}=b$ and $k=1$, $2$ is
spectrally stable. On the other hand, for
$0<d< \frac{a^{2}}{4{\alpha}^{2}}$, the  function $N_{3}(q)$
considered on the nonnegative semi-axis has the minimal value at
\[
q_{0,0}=\frac{a^{2}-4d{\alpha}^{2}}{12d}>0.
\]
Thus $N_{3}(q)$ is monotonically decreasing on the interval $(0, q_{0,0})$
from its zero value at the origin and therefore has the intervals of
negativity, which yields the spectral instability for problem (\ref{spprk}) with
$\frac{a^{2}}{4}=b$ for the small 
values of the diffusion coefficient
$0<d<\frac{a^{2}}{4{\alpha}^{2}}$.
$\hfill\lanbox$

\medskip

In order to  study  spectral stability for   stationary solutions of our integro-differential equation in the case of the Gaussian kernel, we will need the following elementary result.

\smallskip

\noindent {\bf Lemma 8.}
{\it  Assume $\frac{a^2}{4}>b$. Then the  function
$g(s)=-s\ln s+s-\frac{b}{c_{1}^{2}}$, $s\in (0, 1)$ is strictly increasing and it has a unique zero  $s_{0}=s_{0}(a,b)$.}

{\it Proof.} Near the endpoints of the interval $(0,  1)$, 
$g(0^{+})=-\frac{b}{c_{1}^{2}}<0$ and
$g(1^{-})=1-\frac{b}{c_{1}^{2}}>0$ via Lemma 6. On this interval
the function is continuously differentiable such that
$g'(s)=-\ln s>0$, $s\in (0, 1)$. Therefore, it is strictly increasing
and has a unique zero at some point $s_{0}=s_{0}(a,b)\in (0,  1)$.
$\hfill\lanbox$

\medskip

Now we come to prove Theorem~4.

{\it Proof of Theorem 4.} When the kernel is given by (\ref{fi3}), we obtain
\[
\Phi_{k}(p)=\Phi_k(|p|)=d|p|^{2}+c_{k}^{2}e^{-\frac{|p|^{2}}{4 \alpha}}-b, \qquad p\in
{\mathbb R}^{n},  k=1,2.
\]

First we turn our attention to the situation where $\frac{a^{2}}{ 4}>b$. 
In the case of $k=1$, $\Phi_{k}(0)=c_{1}^{2}-b>0$ via Lemma 6 and
we compute
\begin{equation}
\label{der}
\frac{\partial \Phi_{1}(|p|)}{\partial |p|}=2|p|\left(d-
\frac{c_{1}^{2}}{4\alpha}e^{-\frac{|p|^{2}}{4\alpha}}\right), \quad
|p| \in [0, \infty).
\end{equation}
Thus when $d\geq \frac{c_{1}^{2}}{4\alpha}$, the derivative above
is nonnegative, $\Phi_{1}(p)>0$, $p\in {\mathbb R}^{n}$ and hence problem
(\ref{spprk}) with $k=1$ is spectrally stable. On the other hand, when
$0<d<\frac{c_{1}^{2}}{4\alpha}$,  formula (\ref{der}) yields
the single critical point outside the origin for $\Phi_{1}(|p|)$,
\begin{equation}
\label{p0}
|p_{0}|=2\sqrt{\alpha}\sqrt{\ln\frac{c_{1}^{2}}{4\alpha d}}.
\end{equation}
We compute the second derivative
$\frac{\partial^{2}\Phi_{1}}{\partial |p|^{2}}(|p_{0}|)=
\frac{c_{1}^{2}|p_{0}|^{2}}{4 \alpha^{2}}e^{-\frac{|p_{0}|^{2}}{ 4\alpha}}>0$
such that $|p_{0}|$ given by (\ref{p0}) is the minimal point of the
function $\Phi_{1}(|p|)$ on the semi-axis $(0, \infty)$. A straightforward
computation yields $\Phi_{1}(|p_{0}|)=c_{1}^{2}g(s)$ with $g(s)$ defined
in Lemma 8 and $s=\frac{4\alpha d}{c_{1}^{2}}$
such that $s\in (0, 1)$. Thus by means of Lemma 8 when
$d\geq s_{0}(a,b)\frac{c_{1}^{2}}{4\alpha}$ and therefore,
when $s\geq s_{0}(a,b)$,  we have $\Phi_{1}(|p_{0}|)>0$  and problem (\ref{spprk})
with $\frac{a^{2}}{4}>b $ and $k=1$ is spectrally stable.
But when
$0<d< s_{0}(a,b)\frac{c_{1}^{2}}{4\alpha}$ such that
$0<s<s_{0}(a,b)$ problem (\ref{spprk}) with $k=1$ is spectrally unstable
due to the existence of the regions of negativity of $\Phi_{1}(p)$, which
completes the proof of Theorem 4 in the case of $\frac{a^2}{4}>b$.

Now, we study the degenerate case of $\frac{a^{2}}{4}=b$.
Then
\[
\Phi_k(p)=\Phi_k(|p|)=d|p|^{2}+c_{1}^{2}e^{-\frac{p^{2}}{4 \alpha}}-b, \quad p\in
{\mathbb R}^{n},  k=1,2,
\]
and $c_1=c_2=\frac{a}{2}=\sqrt{b}$. Note that $\Phi_{1}(0)=c_{1}^{2}-b=0$ via Lemma 6 and the
derivative with respect to the radial variable is given by (\ref{der}).
Hence when  $d\geq \frac{b}{4\alpha}$ the derivative (\ref{der}) is positive outside the origin such that the function $\Phi_{1}(p)>0$ for
$|p| \in (0, \infty)$, which implies the spectral stability. On the other
hand, when $0<d<\frac{b}{4\alpha}$ the derivative given by (\ref{der})
is negative for small values of $|p|\in (0, \infty)$. Therefore,
problem (\ref{spprk}) with $\frac{a^{2}}{4}=b$ and $k=1$, $2$
is spectrally unstable in this case due to the existence of the regions of
negativity of $\Phi_{1}(p)$.
$\hfill\lanbox$

\medskip

Finally, we prove Theorem 5 in an analogous way as for  Theorem 3.

\smallskip

{\it Proof of Theorem 5.} In the case of the kernel function given by
(\ref{fi4}),
\[
\Phi_{k}(p)=d|p|^{2}+c_{k}^{2}\frac{{\alpha}^{4}}{ ({\alpha}^{2}+|p|^{2})^{2}}-b,
\quad p\in {\mathbb R}^{3}, \quad k=1,2.
\]
By comparing the expression above with formula (\ref{Fik3}) in the case of
$p_{1}^{2}=p_{2}^{2}$ we draw the conclusion that when
$\frac{a^{2}}{4}>b$ and
$d\geq \frac{2x^{*}}{{\alpha}^{2}}$ with $x^{*}$
defined in Lemma 7, the 
function $\Phi_{1}(p)$ is nonnegative in the space of three dimensions
and problem (\ref{spprk}) with $k=1$ is spectrally stable. Moreover, 
$0<d< \frac{2x^{*}}{{\alpha}^{2}}$ yields  regions of
negativity for $\Phi_{1}(p)$ and therefore, spectral instability.
In the degenerate case of $\frac{a^{2}}{4}=b$ the expression
for $\Phi_{k}(p)$, $k=1$, $2$ is given by
\[
d|p|^{2}+\frac{a^{2}}{4}\frac{{\alpha}^{4}}{ ({\alpha}^{2}+|p|^{2})^{2}}-
\frac{a^{2}}{4}, \quad p\in {\mathbb R}^{3}.
\]
Let us compare it with formula (\ref{Fi123}) when $p_{1}^{2}=p_{2}^{2}$.
It follows that when
$d\geq \frac{a^{2}}{2{\alpha}^{2}}$ there is spectral stability
and when $0<d< \frac{a^{2}}{2{\alpha}^{2}}$ the problem
(\ref{spprk}) with $k=1$, $2$ and  $\frac{a^{2}}{4}=b$ is
spectrally unstable. $\hfill\lanbox$


\section*{Acknowledgement} Y.C.\ was partially supported by the NSERC of Canada
(Nos.\ RGPIN-2019-05892, RGPIN-2024-05593).


\end{document}